\documentclass[12pt]{amsart}
   \newtheorem{lemma}{Lemma}[section]
   \newtheorem{theorem}[lemma]{Theorem}
   \newtheorem{remark}[lemma]{Remark}

   \newtheorem{definition}[lemma]{Definition}

\renewcommand{\phi}{\varphi}

\title[Stability of   Stochastic    PDEs]
{Exponential stability of non-autonomous  stochastic   partial
differential equations with finite memory }

\thanks{This work was supported in part by the Natural Science
Foundation of China (No.10171059) and by the NSF Grant 0620539.}

\author{Li Wan }
\address[L.~Wan ]
{1. Department of Mathematics\\ Huazhong University of Science and
Technology \\Wuhan 430074, China \\
2.  Department of Mathematics and Physics\\
Wuhan University of Science and Engineering \\
Wuhan 430073, China } \email[L.~Wan]{wanlinju@yahoo.com.cn}

\author{Jinqiao  Duan }
\address[J.~Duan ]
{Department of Applied Mathematics
\\
Illinois Institute of Technology
Chicago  \\
IL 60616, USA} \email[J.~Duan ]{duan@iit.edu}

\date{August 9, 2007}



\begin{document}

\begin{abstract}

The exponential stability, in both mean square and almost sure
senses, for energy solutions to a
 nonlinear and non-autonomous
stochastic   PDEs with finite memory   is investigated. Various
criteria for stability are obtained. An example is presented to
demonstrate the main results.
\end{abstract}

\maketitle

\section{Introduction}
Recently stochastic partial differential equations  have attracted
a lot of attention, and various results on the existence,
uniqueness and the asymptotic behaviors of the solutions have been
established; see, for example, \cite{CLT,I1,I2,K,LMR,L2,L4,M,WD}.
In particular, stability of solutions has been studied by the
methods of coercivity conditions, the Lyapunov functionals, and
energy
estimates (see \cite{CL,L3,T1}).\\

A few authors have studied  stochastic   partial differential
equations in which the forcing term contains some hereditary
features; see, for example, \cite{CGR,CLT,T3}. These situations
may appear, for instance, when   controlling a system by applying
a force which takes into account not only the present state of the
system but also the history of the solutions.
 The exponential
stability of the mild solutions to the semilinear stochastic delay
evolution equations was discussed by  using Lyapunov functionals
(see \cite{L1,T2}). When   discussing the asymptotic behavior of
solutions, the method by Lyapunov functionals is powerful. However
it is well known that the construction of Lyapunov functionals is
more difficult for functional differential equations such as
differential
equations with memory. \\

The purpose of this paper is to discuss    the mean square
exponential stability
 and almost sure exponential stability of the energy
solutions to the following nonlinear and non-autonomous stochastic
partial differential equation with finite memory:
\begin{eqnarray}
d X(t)&=&[A(t,X(t))+F (t,X_t) ] dt+ G (t,X_t) d W(t),\; t\geq
 0, \end{eqnarray}
\begin{eqnarray*} X(s)=\phi(s)\in L^2(\Omega, C([-r,0],H)),\;\;s\in[-r, 0],\end{eqnarray*}
in which $C:=C([-r,0],H))$ denotes the space of all continuous
functions from $[-r, 0]$ into $H,$ $\phi$ is
$\mathbb{F}_0$-measurable and $A : [0,\infty)\times V\rightarrow
V^*$ and $F :[0,\infty)\times C\rightarrow V^* $ and $G
:[0,\infty)\times C\rightarrow L^0_Q(K,H)$
 are continuous.  \\

The contents of this paper are as follows. In Section 2 we present
some preliminaries and consider the existence of energy solutions
(see Definition 2.1).  In Section 3 we consider stability of the
nonlinear and non-autonomous stochastic  partial differential
equations with finite memory. In Section 4 we present an example
which illustrates the main results in this paper.

\section{Preliminaries}

Let $V$, $H$ and $K$ be separable Hilbert spaces  and let $L(K,H)$
be the space of all bounded linear operator from $K$ to $H$. We
denote the norms of elements in $ V, H, K$ and $L(V,H)$ by
 $\|\cdot\|, |\cdot|_2, |\cdot|_K$ and
$|\cdot|$ respectively. And $|\cdot|_*$ denotes the norm of $
V^*$, $<\cdot,\cdot>$ denotes the duality between $V$ and $V^*.$
    Suppose that $V$ and $H$ satisfy
$$V\subset H\equiv H^*\subset V^*,$$ where $V$ is a dense subspace
of $H$ and the injections are continuous with
\begin{eqnarray}\lambda_1|v|^2_2\leq\| v\|^2, \;\;\lambda_1>0,\;\;v\in V. \end{eqnarray}

  We are given a
$Q$-Wiener process in the complete probability space
$(\Omega,\mathcal{F},P,\{\mathcal{F}_t\}_{t\geq 0})$ and have
values in $K$, i.e. $W(t)$ is defined as
\begin{equation*}W(t)=\sum^{\infty}_{n=1}
\sqrt{\lambda_n}B_n(t)e_n, t\geq 0,
\end{equation*}where
$B_n(t),(n=1,2,\cdots)$ is a sequence of  real value standard
Brownian motions mutually independent on
$(\Omega,\mathcal{F},P,\{\mathcal{F}_t\}_{t\geq 0}), \lambda_n\geq
0,(n=1,2,\cdots)$ are nonnegative real numbers such that
$\sum_{n\geq 1}\lambda_n< \infty, \{e_n\}_{n\geq 1}$ is a complete
orthonormal  basis in $K$, and $Q\in \mathcal{L}(K,K)$ is the
incremental covariance operator of the process $W(t)$, which is a
symmetric nonnegative trace class operator defined as
\begin{equation*} Qe_n=\lambda_n e_n, n=1,2,\cdots
\end{equation*} Let  $L^0_Q(K,H)$ be the space of all bounded linear operators from
$K$ to $H$ with the following condition:
$$\|\xi\|^2_{L^0_Q}=tr(\xi Q\xi^*)<\infty, \;\;\xi\in L(K,H).$$
Let $M^2(-r,T; V)$ denote the space of all $V$ -valued measurable
functions defined on  $[-r,T]\times \Omega$ with
$E\int^T_{-r}\|X(t)\|^2dt<\infty.$\\

First we give the definition of an energy solution to (1).

\begin{definition} Stochastic process X(t) on $(\Omega,\mathcal{F},P,\{\mathcal{F}_t\}_{t\geq 0})
$ is called an energy solution to (1) if the following conditions
are satisfied:
\begin{itemize}
\item [(i)] $X(t)\in M^2(-r,T; V)\cap L^2(\Omega; C(-r,T; H)),
T>0,$

 \item [(ii)]  the following equation holds in $V^*$  almost surely, for $t\in [0,T
 ),$ \begin{eqnarray*}   X(t)&=&X(0)+\int^t_0[A(s,X(s))+F(s,X_s) ] ds+ \int^t_0G(s,X_s) d W(s),\; t\geq
 0, \end{eqnarray*}
\begin{eqnarray*} X(s)=\phi(s) ,\;\;s\in[-r, 0],\end{eqnarray*}

 \item [(iii)]  the following stochastic energy equality holds:
\begin{eqnarray}   |X(t)|^2_2&=& |X(0)|^2_2+2\int^t_0 <X(s), A(s,X(s))+F(s,X_s) >ds
\nonumber\\&&+\int^t_0 \|G(s,X_s)\|^2_{L^0_Q}ds+
2\int^t_0<X(s),G(s,X_s) d W(s)>.
\end{eqnarray}
\end{itemize}
\end{definition}

 In order to guarantee the existence and
uniqueness of energy solution $X(t)$ to (4), we need the following
 conditions:

\begin{itemize}
\item [(A1)] (Monotonicity and coercivity) There exist $\alpha>0$
and $\lambda\in \mathbb{R}$ such that for a.e. $ t\in (0,T ),$
$$-2<A(t,u)-A(t,v), u-v>+ \lambda|u-v|^2_2 \geq
\alpha\|u-v\|^2.$$
 \item [(A2)]  (Measurability) For any $v \in V,$ the mapping $t \in (0,T )\rightarrow A(t, v) \in V^*$ is measurable.
 \item [(A3)](Hemicontinuity) The next mapping is continuous for any $u, v,w \in V,$ a.e. $t \in (0,T
 ):$ $$\mu\in \mathbb{R}\rightarrow <A(t, u+\mu v), w>\in
 \mathbb{R}.$$
 \item [(A4)](Boundedness) There exists $c>0$ such that for any $  v  \in V,$ a.e. $t \in (0,T
 ):$
$$|A(t,v)|_*\leq c\|v\|.$$
 \item [(A5)] (Lipschitz condition) There exists  $c_1 > 0 $ such that for any
 $\xi, \eta\in C$ and $F (t,0)\in L^2([-r,0]\times\Omega, V^*)$
 and $G (t,0)\in L^2([-r,0]\times\Omega,H),$
 $$\|F (t,\xi)-F (t,\eta)\|_*\leq c_1 |\xi-\eta|_C,\;\; \|G (t,\xi)-G (t,\eta)\|_ {L^0_Q}\leq c_1
 |\xi-\eta|_C.$$
\end{itemize}

We have the following result on the energy  solution to (4) (see
\cite{CGR,L4,P}):

\begin{theorem}Suppose that conditions $(A1)-(A5)$ are satisfied. Then there exists
a unique energy solution $X(t)$ to (4). Furthermore, the following
identity holds:
\begin{eqnarray*}   \frac{d}{dt}E|X(t)|^2_2&=&  2 E<X(s), A(t,X(t))+F (t,X_t)
>
+E\|G(t,X_t)\|^2_{L^0_Q}, t\geq t_0.
\end{eqnarray*}
\end{theorem}

Throughout this paper, we assume the existence of the energy
solutions to (1) with $E\|\phi\|^2_{C}<\infty.$

\section{Main results}

In this section we discuss the stochastic   partial differential
equations with finite memory which are the examples of (1).  Let
$\rho, \tau : [0,\infty) \rightarrow [0, r]$ be continuous
functions  and $r>0$ is a constant. Assume that
$A:[0,\infty)\times V\rightarrow V^*$ and $f  :[0,\infty)
\rightarrow V^*$ are Lebesgue measurable. Let $g:[0,\infty)\times
H\rightarrow H$ and $h:[0,\infty)\times H\rightarrow L_Q(K,H)$ be
Lipschitz continuous uniformly in t.\\

We consider the following stochastic partial differential equation
with finite memory:

\begin{eqnarray}
 d X(t)&=&[A(t,X(t))+f (t) ] dt+g(t,X(t-\rho(t))dt\nonumber\\&&+h(t,X(t-\tau(t))dW(t),\; t\geq
 0,
\end{eqnarray}
 with the initial condition \begin{eqnarray*} X(s)=\phi(s)\in L^2(\Omega, C([-r,0],H)),\;\;s\in[-r, 0],\end{eqnarray*}
where $f \in L^2([0,\infty),V^*).$ Set $F_1(t,
\psi)=g(t,\psi(-\rho(t)))$ and $G(t, \psi)=h(t,\psi(-\tau(t)))$
for any $\psi \in C.$  Then   (4) can be viewed as a stochastic
functional partial differential equation (1) with $F(t,
\psi)=F_1(t, \psi)+f (t).$ Note that we do not require that
$\rho(t)$ and $\tau(t)$ are
differentiable functions.\\

Now we first state the  following result on  exponential stability
in mean square.

\begin{theorem}Suppose that  Eq.  (4) satisfies the following conditions:
\begin{itemize}
\item [(B1)]conditions $(A2)-(A4)$ holds and there exist
$\delta_1>0$ and a continuous, integrable function $\alpha_1(t)>0$
such that for a.e. $t\in [0,\infty), u, v \in V,$
$$-2<A(t,u)-A(t,v), u-v>+ \alpha_1(t)|u-v|^2_2\geq
\delta_1\|u-v\|^2;$$
 \item [(B2)]  there exist integrable functions $\alpha_2,
 \beta_2: [0,\infty) \rightarrow \mathbb{R}^+ $ such that
 $$ |g(t,u)|^2_2\leq (\delta_2 +\alpha_2(t)) |u |^2_2+ \beta_2(t), for\;\; \delta_2\geq 0, u\in H;  $$
 \item [(B3)]   there exist integrable functions $\alpha_3,
 \beta_3: [0,\infty) \rightarrow \mathbb{R}^+ $ such that
$$\|h(t,u)\|^2_{L^0_Q}\leq (\delta_3 +\alpha_3(t)) |u |^2_2+ \beta_3(t), for\;\; \delta_3\geq 0, u\in H;  $$
 \item [(B4)]  there exists $ \sigma_1  > 0$ such that
 $$\int^\infty_0 e^{\sigma_1 t} |f_\beta(t)|^2_*   dt<\infty,\;\;\int^\infty_0 e^{\sigma_1 t}
 \beta_i(t)dt<\infty,\;i=2,3;$$
  \item [(B5)] $\delta_1\lambda_1 >
2\sqrt{\delta_2}  +   \delta_3,$ where $\lambda_1$ is defined by
(2).
\end{itemize} Then for any energy solution X(t) to (4), there
exist $\sigma\in (0, \sigma_1)$ and $B> 1$ such that
\begin{eqnarray}E|X(t)|^2_2\leq Be^{-\sigma t}  ,  t \geq 0.\end{eqnarray} In other words the
energy solution X(t) to (4) converges to zero exponentially in
mean square as $t\rightarrow\infty$.

\end{theorem}
\begin{proof}From $(B5)$, there exists  $\gamma_1>0$
such that
 $$( \delta_1-\gamma_1)\lambda_1 >
2\sqrt{\delta_2}  +   \delta_3 .$$ Thus we can take $\gamma_2>0$
such that  $$( \delta_1-\gamma_1)\lambda_1 >\gamma_2+
   \frac{\delta_2  }{\gamma_2} +
  \delta_3>
2\sqrt{\delta_2}  +   \delta_3.$$ Furthermore,  we can choose
$\sigma\in (0, \sigma_1)$
  such that
   $$( \delta_1-\gamma_1)\lambda_1 > \sigma+\gamma_2+
  e^{\sigma r}  \frac{\delta_2  }{\gamma_2} + e^{\sigma r}
  \delta_3.$$

 Now define the function $f_1: [0,\infty)\times V\rightarrow
V^*$ by
\begin{eqnarray*}
  f_1(t, v)=A(t,v)+f (t) , v\in V,\;t\geq 0.
\end{eqnarray*}  From $(B1)$ and $f  \in
L^2([0,\infty),V^*),$ it follows that
\begin{eqnarray*}
 2<f_1(t, v),v> &\leq&
 \alpha_1(t)| v|^2_2-\delta_1\| v\|^2+2<f (t), v>\\&\leq&
 \alpha_1(t)| v|^2_2-\delta_1\| v\|^2+\gamma_1\|v\|^2 +\gamma_1^{-1}|f (t)|^2_*\\&\leq&
 [(-\delta_1+\gamma_1)\lambda_1+\alpha_1(t)]| v|^2_2 +\gamma_1^{-1}|f (t)|^2_*\\&=&
 [ -a+\alpha_1(t)]| v|^2_2 +\beta_1(t),\;\; v\in V,
\end{eqnarray*}where $a=( \delta_1-\gamma_1)\lambda_1, \beta_1(t)=\gamma_1^{-1}|f (t)|^2_*.$

Set
$$\theta(t )=\alpha_1(t)  +
  e^{\lambda r}  \frac{ \alpha_2(t )}{\gamma_2}
      +e^{\lambda r}  \alpha_3(t  ),\;\;\beta(t)=\beta_1(t)+\frac{\beta_2(t)}{\gamma_2}+\beta_3(t).$$
It follows   from $(B1)-(B4)$ and $f \in L^2([0,\infty),V^*) $
that
\begin{eqnarray*}R_1=\int^
\infty_0 \theta(s)  ds <\infty , R_2=\int^ \infty_0 \beta(s) ds
\leq  R_3=\int^ \infty_0 e^{\sigma_1 s}\beta(s) ds<\infty.
\end{eqnarray*}

 Set\begin{eqnarray}
     K(t)=
    \bigg\{\begin{array}{cc}
    E|X(t)|^2_2e^{\sigma t}\exp\bigg(-\int^t_0[\theta(s)+e^{\sigma s}\beta(s)]ds\bigg)
    ,& t \geq 0,\\E|X(t)|^2_2e^{\sigma t} ,&-r\leq t
  < 0.
\end{array}
\end{eqnarray} It is clear that $K(t)$ is continuous on $[-r, \infty) $ and
\begin{eqnarray} \frac{ d K(t)}{dt}
 \nonumber &=& e^{\sigma t}\exp\bigg(-\int^t_0[\theta(s)+e^{\sigma s}\beta(s)]ds\bigg)\bigg\{ \sigma E|X(t)|^2_2
 -[\theta(t)+e^{\sigma t}\beta(t)]
  E|X(t)|^2_2\nonumber\\&&+2E<f_1(t, X(t)), X(t)>+
  2E<g(t, X(t-\rho(t))), X(t) >\nonumber\\&&+ E\|h(t,X(t-\tau(t)))\|^2_{L^0_Q}\bigg\}.\end{eqnarray}
  From (B2)
and (B3), it follows that  for $t\geq 0$
\begin{eqnarray}
\frac{ d K(t)}{dt}
    \nonumber&\leq & e^{\sigma t}\exp\bigg(-\int^t_0[\theta(s)+e^{\sigma s}\beta(s)]ds\bigg)\bigg\{ \sigma E|X(t)|^2_2
-[\theta(t)+e^{\sigma t}\beta(t)]E|X(t)|^2_2
    \nonumber\\&&+[ -a+\alpha_1(t)]E| X(t) |^2_2
 +\beta_1(t)+\gamma_2 E|X(t)|^2_2+
   \frac{ E|g(t,X(t-\rho(t)))|^2_2}{\gamma_2}     \nonumber\\&&   + E\|h(t,X(t-\tau(t)))
   \|^2_{L^0_Q}\bigg\}
      \nonumber  \\ &\leq&   [ -a+\alpha_1(t)+\sigma+\gamma_2-\theta(t) ] K(t )
   + e^{\sigma t}
    \beta(t) -   e^{ \sigma t}    \beta(t)
    K(t )  \nonumber \\&&     +e^{\sigma t}
   \exp\bigg(-\int^t_0[\theta(s)+e^{\sigma s}\beta(s)]ds\bigg)
   \frac{\delta_2 +\alpha_2(t)}{\gamma_2}
   E|X(t-\rho(t))|^2_2\nonumber\\&&
   + e^{\sigma t}
   \exp\bigg(-\int^t_0[\theta(s)+e^{\sigma s}\beta(s)]ds\bigg) (\delta_3 +\alpha_3(t)) E|X(t-\tau(t))|^2_2
   .
 \end{eqnarray}

Now we claim  that \begin{eqnarray}K(t)\leq
 1+\sup_{t\in[-r,0]}\{E|X(t)|^2_2\}=M \;\;\;\mbox{ for all } t\geq
0.
\end{eqnarray} In fact, if
(9) is false, then there exists $t_1>0$ such that for $\forall
\varepsilon>0$,
\begin{eqnarray}K(t)\leq M, \;\; 0\leq t<t_1 ,\;\;K(t_1)=M;\;\;K(t)>M, \;\; t_1<t<t_1+\varepsilon.\end{eqnarray}
Since $\frac{ d K(t)}{dt}$ exists (by (7)), we see that
\begin{eqnarray}\frac{d}{dt}K(t_1) \geq 0.
\end{eqnarray} From
(8), it follows that
\begin{eqnarray}
 \frac{ d }{dt}K(t_1)
  \nonumber  &\leq &      [ -a+\alpha_1(t_1)+\sigma+\gamma_2 -\theta(t_1)
    ]K(t_1)
 \nonumber\\&&  +
  e^{\sigma t_1}\exp\bigg(-\int^{t_1}_0[\theta(s)+e^{\sigma s}\beta(s)]ds\bigg) \frac{\delta_2 +\alpha_2(t_1)}{\gamma_2}
   E|X(t_1-\rho(t_1))|^2_2\nonumber\\&&
   +e^{\sigma t_1}\exp\bigg(-\int^{t_1}_0[\theta(s)+e^{\sigma s}\beta(s)]ds\bigg)
    (\delta_3 +\alpha_3(t_1)) E|X(t_1-\tau(t_1))|^2_2      .
 \end{eqnarray}

 We consider the following three different cases:\\

i)  Suppose that $t_1-\rho(t_1)\geq 0, t_1-\tau(t_1)\geq 0.$ From
(10) and  (12), it follows that
\begin{eqnarray*}
\frac{ d }{dt}K(t_1)
 &\leq  &     [ -a+\alpha_1(t_1)+\sigma+\gamma_2 -\theta(t_1) ]
 K(t_1)
 \nonumber\\&&  +
  e^{\sigma \rho(t_1)}\exp\bigg(-\int^{t_1}_{t_1-\rho(t_1)}[\theta(s)+e^{\sigma s}\beta(s)]ds\bigg) \frac{\delta_2 +\alpha_2(t_1)}{\gamma_2}
   K(t_1-\rho(t_1)) \nonumber\\&&
   +e^{\sigma \tau(t_1)}\exp\bigg(-\int^{t_1}_{t_1-\tau(t_1)}[\theta(s)+e^{\sigma s}\beta(s)]ds\bigg) (\delta_3 +\alpha_3(t_1))
 K(t_1-\tau(t_1))
   \end{eqnarray*}\begin{eqnarray*} \\&\leq &    [ -a+\alpha_1(t_1)+\sigma+\gamma_2 -\theta(t_1) ]
 M
 +
  e^{\sigma r}  \frac{\delta_2 +\alpha_2(t_1)}{\gamma_2}
  M
   \nonumber\\&&  +e^{\sigma r} (\delta_3 +\alpha_3(t_1))
M      \\&\leq &     [ -a+ \sigma+\gamma_2+
  e^{\sigma r}  \frac{\delta_2  }{\gamma_2} + e^{\sigma r}  \delta_3 ]
 M<0 .
 \end{eqnarray*}

ii)  Suppose that $t_1-\rho(t_1)<0, t_1-\tau(t_1)\geq 0.$ Then
$t_1-\rho(t_1)>-r .$ From (10) and  (12),  it follows that
\begin{eqnarray*}
\frac{ d }{dt}K(t_1)
 &\leq &     [ -a+\alpha_1(t_1)+\sigma+\gamma_2 -\theta(t_1) ]
 K(t_1)
 \nonumber\\&&  +
  e^{\sigma \rho(t_1)}\exp\bigg(-\int^{t_1}_{0}[\theta(s)+e^{\sigma s}\beta(s)]ds\bigg) \frac{\delta_2 +\alpha_2(t_1)}{\gamma_2}
   K(t_1-\rho(t_1)) \nonumber\\&&
   +e^{\sigma \tau(t_1)}\exp\bigg(-\int^{t_1}_{t_1-\tau(t_1)}[\theta(s)+e^{\sigma s}\beta(s)]ds\bigg) (\delta_3 +\alpha_3(t_1))
 K(t_1-\tau(t_1))
     \\&\leq  &    [ -a+ \sigma+\gamma_2+
  e^{\sigma r}  \frac{\delta_2  }{\gamma_2} + e^{\sigma r}  \delta_3 ]
 M<0 .
 \end{eqnarray*}
iii)  Suppose that $t_1-\rho(t_1)<0, t_1-\tau(t_1)< 0.$ Then
$t_1-\rho(t_1)>-r, t_1-\tau(t_1)>-r .$ Similarly, it follows that
\begin{eqnarray*}
\frac{ d }{dt}K(t_1)
 &\leq &       [ -a+ \sigma+\gamma_2+
  e^{\sigma r}  \frac{\delta_2  }{\gamma_2} + e^{\sigma r}  \delta_3 ]
 M<0 .
 \end{eqnarray*}Thus,  one  obtains $$\frac{d}{dt}K(t_1)< 0 ,$$
 which
  contradicts with (11). Hence, (9) holds true and we obtain
$$E|X(t)|^2_2\leq e^{-\sigma
t}\exp\bigg( \int^t_0[\theta(s)+e^{\sigma s}\beta(s)]ds\bigg) M
\leq e^{-\sigma t}B ,  t \geq 0,$$ where $B=e^{R_1+R_3} M.$
 The proof is complete.
\end{proof}
 \begin{remark} Comparing  with \cite{T4}, our method does not assume that $\rho$ and $\tau$ are
 differentiable functions. Hence, our method can be applied to more general
 stochastic
 partial differential equations with memory.

\end{remark}

 Next, we state the result of almost sure exponential stability.
\begin{theorem}  Suppose that all the conditions of Theorem 3.1 are
satisfied. If the following additional condition is satisfied:
\begin{itemize}
   \item [(B6)]  $\alpha_i(t)$ and $e^{\sigma t} \beta_i(t) (i=1,2,3) $are bounded functions,
   where $\beta_1(t)=\gamma_1^{-1}|f (t)|^2_*.$
\end{itemize}
 Then
   there exists
$T(\omega)>0$ such that for all $t>T(\omega)$
\begin{eqnarray*} |X(t)|^2_2\leq e^{ \sigma /2 }e^{-\sigma t/2}
\end{eqnarray*}
 with probability one.

\end{theorem}
\begin{proof}   Let $N_1$ and $N_2$ be positive integers such
that $$N_1-\rho(N_1)\geq N_1-r\geq 1 ,\;\;N_2-\tau(N_2)\geq
N_2-r\geq 1 .$$ Let $N>N_3=\max(N_1, N_2) $ and $I_N=[N,N+1].$

Set
$$\alpha (t)= \alpha_1(t)+ \gamma_2
 +
 \frac{\delta_2 +\alpha_2(t)}{\gamma_2} e^{ \sigma r}+32(\delta_3 +\alpha_3(t))e^{ \sigma r}  , $$
  $$\beta (t)=
2\beta_1(t)+ \frac{2\beta_2(t)}{\gamma_2}+ 64\beta_3(t).$$ It
follows from (B6) that there exists $B_1>0$ such that $$\alpha
(t)+e^{\sigma t}\beta (t)\leq B_1.$$

Then
 we obtain from (3)
\begin{eqnarray*}
 {E}  \sup_{t\in I_N}|X(t)|^2_2&\leq & E|X(N)|^2_2+2{E}  \sup_{t\in
 I_N}\int^t_N
<X(s), A(s,X(s))+ f (s) >ds \nonumber\\&&+2{E}  \sup_{t\in
 I_N}\int^t_N
<X(s), g(t,X(s-\rho(s))  >ds \nonumber\\&&+ {E}  \sup_{t\in
 I_N}\int^t_N
\|h(t,X(s-\tau(s))\|^2_{L^0_Q}ds\nonumber\\&&+ 2{E}  \sup_{t\in
 I_N}\int^t_N<X(s),h(t,X(s-\tau(s))) d
W(s)>.\end{eqnarray*} Then, by the   Burkholder-Davis-Gundy
inequality (\cite{Mao,RY}),
\begin{eqnarray*}
 &&2{E}  \sup_{t\in
 I_N}\int^t_N<X(s),h(t,X(s-\tau(s))) d
W(s)>\\ &\leq & \frac{1}{2}{E}  \sup_{t\in I_N}|X(t)|^2_2+32
 \int^{N+1}_N
E\|h(t,X(s-\tau(s))\|_{L^0_Q}ds.\end{eqnarray*} Therefore it
follows from $(B1)-(B3)$ and (5) that
\begin{eqnarray*}
 {E}  \sup_{t\in I_N}|X(t)|^2_2&\leq & 2E|X(N)|^2_2+4{E}  \sup_{t\in
 I_N}\int^t_N
 <X(s), A(s,X(s))+ f (s) >  ds \nonumber\\&&+4{E}
\sup_{t\in
 I_N}\int^t_N
<X(s), g(t,X(s-\rho(s))  >ds \nonumber\\&&+ 64    \int^{N+1}_N
E\|h(t,X(s-\tau(s))\|^2_{L^0_Q}ds
\\&\leq & 2E|X(N)|^2_2+2
 \int^{N+1}_N  \alpha_1(s) E| X(s)|^2_2 +\beta_1(s)
ds \nonumber  \\&&+2 \int^{N+1}_N \gamma_2 E|X(s)|^2_2+
   \frac{\delta_2 +\alpha_2(s)}{\gamma_2}  E|X(s-\rho(s))|^2_2+
\frac{\beta_2(s)}{\gamma_2}  ds
\nonumber\end{eqnarray*}\begin{eqnarray*}\\&&+ 64  \int^{N+1}_N
(\delta_3 +\alpha_3(s)) E|X(s-\tau(s))|^2_2+ \beta_3(s)ds\\&\leq &
 2B e^{-\sigma N} +
 \int^{N+1}_N 2B e^{-\sigma s}\bigg\{  \alpha (s)  +\beta (s)e^{ \sigma
s}\bigg\}ds \leq    2B  e^{-\sigma N}(1 +\frac{B_1}{\sigma})
 .\end{eqnarray*}Let $\varepsilon_N $ be any fixed positive real number. Then,
\begin{eqnarray*}  \mathbb{P}\{ \sup_{t\in I_N}|X(t)|^2 _2
>  \varepsilon_N^2 \}  &\leq& \frac{{E}  \sup_{t\in I_N}|X(t)|^2_2}{\varepsilon_N^2}
 \leq \frac{2B  e^{-\sigma N}(1 +\frac{B_1}{\sigma}) }{\varepsilon_N^2}.
\end{eqnarray*}

For each integer $N $, choosing $\varepsilon_N^2 =  e^{-\frac{
\sigma N}{2}},$ then
\begin{eqnarray*}\mathbb{P}\{ \sup_{t\in I_N}|X(t)|^2_2
>   e^{-\frac{
\sigma N}{2}} \}  &\leq&   2B  e^{-\sigma N/2}(1
+\frac{B_1}{\sigma})  .
\end{eqnarray*}
From Borel-Cantelli's lemma (\cite{Ash, Mao}), it follows that
there exists $T(\omega)>0$ such that for all $t>T(\omega)$
\begin{eqnarray*} |X(t)|^2_2\leq e^{ \sigma /2 }e^{-\sigma t/2},
a.s.
\end{eqnarray*}
 The proof is complete.

\end{proof}

\section{ An Example}

In this section we present an example which illustrates the main
results. We consider a sufficient condition for any energy
solution to a stochastic   heat equation with finite memory to
converge to zero  exponentially in mean square and
 almost surely exponentially.
Let $A=\gamma_1\frac{\partial^2}{\partial x^2} $,where $\gamma_1
> 0, H = L^2(0, \pi)  $ and  $$H^1 _0 = \{u \in L^2(0, \pi): \frac{\partial u}{\partial x }  \in L^2(0, \pi),
u(0) = u(\pi) = 0\}, $$ $$H^2 = \{u \in L^2(0, \pi):
\frac{\partial u}{\partial x } , \frac{\partial^2}{\partial x^2}
\in L^2(0, \pi), u(0) = u(\pi) = 0\}.$$ Operator $A$ has the
domain $D(A) = H^1 _0 \cap H^2.$ Define the norms of two spaces
$H$ and $H^1 _0$ by $|\xi |^2_2 =  \int^\pi_0 \xi^2(x) dx  $ for
any $\xi\in  H$ and $\|u\|^2= \int^\pi_0  ( \frac{\partial
u}{\partial x }  )^2 dx  $  for any $u \in  H^1 _0 ,$
respectively. Then it is known that
$$<Au, u> \leq -\gamma_1\|u\|^2, u\in H^1 _0 .$$
 Let  $ \rho(t)$ and   $ \tau(t)$ be the non-differentiable function defined by
    $$  \rho(t)= \frac{1}{1+|\sin t|},  \;\; \tau(t)= \frac{1}{1+|\cos t|},  \;\;t\geq  0. $$
  Consider
the following stochastic   heat equation with finite memory
$\rho(t)$ and $\tau(t)$: \begin{eqnarray}dX(t) = AX(t)+g(t, X
(t-\rho(t)))dt +h(t, X (t-\tau(t)))dw(t ),\end{eqnarray} $$ X(t,
0) = X(t,\pi) = 0, t\geq  0, X(s, x) =  \phi(s,x), s \in [-1/2 ,
0], x\in [0,\pi],$$
  $$\phi\in C ([-1/2 , 0],  L^2(0, \pi)) ,  \|\phi\|_C<\infty, $$
   where $w(t)$ is the
one-dimensional standard Wiener process, $$g(t, y )
=(b_1+k_1(t))y+e^{-kt}p,\;\; h(t, y) = (b_2 +k_2(t))y,$$ for any
$y\in H, t \geq 0,$   $p\in H$ with $|p|_2 <\infty, $ $k_1, k_2 :
[0,\infty) \rightarrow\mathbb{R}^+$ are continuous functions,
$b_1,  b_2$ and $k$ are  positive real numbers.  It is easy to
obtain
$$|g(t, y )|^2_2\leq4(b^2_1 +k^2_1(t))|  y
|^2_2+2e^{-2kt}|p|^2_2,$$
$$\|h(t, y )\|^2_{L^0_Q}\leq4(b^2_2 +k^2_2(t))|  y
|^2_2.$$ Note that $\lambda_1=1, \delta_1=2\gamma_1,
\delta_2=4b^2_1, \delta_3=4b^2_2, \alpha_2(t)=4k^2_1(t),
\alpha_3(t)=4k^2_2(t), \beta_2(t)=2e^{-2kt}|p|^2_2, \beta_3(t)=0 $
and taking $\alpha_1(t)=1, \sigma_1=2k.$ Now we suppose that
$$2\gamma_1>4b_1+4b^2_2$$ and $k^2_1(t), k^2_2(t)$ are decreasing, bounded and integrable
functions. Then, it follows from Theorem  3.1 and Theorem  3.3
that any energy solution $X(t)$ to (13) converges to zero
exponentially in mean square and  almost surely exponentially as
$t\rightarrow\infty.$


\begin{thebibliography}{100}
\bibitem{Ash} R. B. Ash,  Probability and Measure Theory.
Academic Press, San Diego, Second Edition, 2000.

\bibitem{CL}  T. Caraballo, K. Liu, On exponential stability criteria of
stochastic partial differential equations, Stochastic Process.
Appl. 83 (1999) 289-301.

\bibitem{CLT}  T. Caraballo, K. Liu, A. Truman,
Stochastic functional partial differential equations; Existence,
uniqueness and asymptotic decay property, Proc. R. Soc. Lond. Ser.
A Math. Phys. Eng. Sci. 456 (2000) 1775-1802.

\bibitem{CGR}  T. Caraballo, M.J. Garrido-Atienza, J. Real, Existence and
uniqueness of solutions for delay stochastic evolution equations,
Stoch. Anal. Appl. 20 (2002) 1225-1256.

\bibitem{CLT}  T. Caraballo, J. Lang, T. Taniguchi, The exponential behaviour
and stabilizability of stochastic 2D-Navier-Stokes equations, J.
Differential Equations (2002) 714-737 .

\bibitem{I1} A. Ichikawa, Stability of semilinear stochastic evolution
equations, J. Math. Anal. Appl. 90 (1982) 12-44.

\bibitem{I2}  A. Ichikawa,
Absolute stability of a stochastic evolution equation, Stochastics
11 (1983) 143-158.

\bibitem{K}  A. Kwiecinska, Stabilization of evolution equation by noise,
Proc. Amer. Math. Soc. 130 (2001) 3067-3074.

\bibitem{LMR} G. Leha, B. Maslowski, G. Ritter, Stability of solutions to
semilinear stochastic evolution equations, Stoch. Anal. Appl. 17
(1999) 1009-1051.

\bibitem{L1}  K. Liu, Lyapunov functionals and asymptotic stability of
stochastic delay evolution equations, Stoch. Stoch. Rep. 63 (1998)
1-26.

\bibitem{L2} R. Liu, V. Mandrekar, Stochastic semilinear evolution
equations: Lyapunov functions, stability, and ultimate
boundedness, J. Math. Anal. Appl. 212 (1997) 537-553.

\bibitem{L3}  K. Liu, X. Mao, Exponential stability of non-linear
stochastic evolution equations, Stochastic Process. Appl. 78
(1998) 173-193.

\bibitem{L4} K. Liu, Stability of infinite dimensional stochastic differential
equations with applications, Chapman \& Hall/CRC, (2006).

\bibitem{M}  B. Maslowski, Stability of semilinear equations with boundary
and pointwise noise, Ann. Sc. Norm. Super. Pisa Cl. Sci. IV (1995)
55¨C93.

\bibitem{Mao} X. Mao, Exponential Stability of Stochastic Differential Equations, Marcel Dekker, New York, 1994.

\bibitem{RY} D. Revuz and M. Yor,
 Continuous Martingales and Brownian Motion.
Springer-Verlag, New York, 1999, Second Edition.

\bibitem{P}  E. Pardoux, ¡äEquations aux d¡äeriv¡äees partielles
stochastiques non lin¡äeaires monotones, Thesis, Universit¡äe
Paris XI,  1975.

\bibitem{T1} T. Taniguchi, Asymptotic stability theorems of semilinear
stochastic evolution equations in Hilbert spaces, Stoch. Stoch.
Rep. 53 (1995) 41-52.

\bibitem{T2}T. Taniguchi, Almost sure exponential stability for
stochastic partial functional differential equations, Stoch. Anal.
Appl. 16 (1998) 965-975.

\bibitem{T3}T. Taniguchi, K. Liu, A. Truman, Existence, uniqueness, and
asymptotic behavior of mild solution to stochastic functional
differential equations in Hilbert spaces, J. Differential
Equations 181 (2002) 72-91.

\bibitem{T4}T. Taniguchi, The exponential stability for stochastic delay partial
differential equations, J. Math. Anal. Appl. 331 (2007) 191-205.

\bibitem{WD}  E. Waymire, J. Duan  (Eds.), Probability and
Partial Differential Equations in Modern Applied Mathematics. IMA
 140,  Springer-Verlag, New York, 2005.

\end{thebibliography}
\end{document}